\documentclass[preprint,authoryear,12pt]{elsarticle}
\usepackage{amssymb}
\usepackage{amsmath}
\usepackage{pgfplots}
\newtheorem{theorem}{Theorem}

\newenvironment{proof}{\noindent{\bf Proof. }}{\hfill\fbox{}\vspace*{3mm}}

\def\al{\alpha}

\def\b{\bf}

\def\bt{\beta}

\def\ga{\gamma}
\def\gath{\gamma\theta}

\def\la{\lambda}

\def\om{\omega}

\def\pa{\partial}

\def\qed{\hfill\fbox{}\vspace*{3mm}}

\def\th{\theta}
\def\dfrac{\displaystyle \frac}

\def\dint{\displaystyle \int}

\def\b1{{\bf{1}}}

\begin{document}
\begin{frontmatter}
\title{Economic production quantity models for imperfect product and service with rework}
\author{Allen H. Tai}

\address{Department of Applied Mathematics, The Hong Kong Polytechnic University, Hung Hom, Hong Kong}

\begin{abstract}
When imperfect quality products are produced in a production process, rework may be performed to make them become serviceable.
In an inventory system, items may deteriorate. 
Selling deteriorated items to customers will create negative impact on corporate image.
In this paper, two economic production quantity (EPQ) models are proposed for deteriorating items
with rework process.
A single production-rework plant system and a system consists of $n$ production plants and one rework plant are considered. 
Approximated analytic results are obtained and numerical examples are provided to illustrate the solution procedure.

\end{abstract}


\end{frontmatter}

\section{Introduction}
Serving high quality products and providing good service can always attract customers and keep them coming back.
However, in reality, production processes are often imperfect.
For economic and environmental reasons, imperfect quality items are reworked to become serviceable again.
Due to unsuitable inventory condition or other reasons, items stored in inventory face deterioration.
In order to provide good service, inspection may be carried out to screen out deteriorated items.
However, when inventory level is huge, full-scale inventory inspection is not possible and only part of deteriorated items can be screen out.
The remaining deteriorated items will then be sold to customers. This will lower customer satisfaction and is harmful to corporate image.

In this paper, we consider two models in which the products produced are imperfect and the service provided to customers is also imperfect.
The products are imperfect in two ways. Firstly, during the production process, imperfect quality items may be produced.
Good quality items are stocked and sold to customers immediately. 
Imperfect quality items are stocked separately and scheduled for rework.
Secondly, items may be deteriorated in the inventory.
Hence, inspection processes are carried out in the inventory.
The service is imperfect in two ways. Firstly, inspections for deteriorated items in inventory are imperfect.
Deteriorated items may be sold to customers. Secondly, shortage is allowed and unsatisfied demands are backlogged.
Since not all customers accept late delivery, partial backlogging is also considered in the paper.

In the first model, we consider a single production plant system. 
The plant is also capable of rework processing.
Part of the imperfect quality items are recovered and ready to be sold to customers.
At the beginning of the production process, the amount of backlog is made up. 
The production process continues until the total number of items produced reaches the economic production quantity $Q^*$.
Then rework process starts until all imperfect quality items are processed.
Demands are then satisfied by the serviceable items in the inventory.
After the inventory level becomes zero, further demands are backlogged and satisfied at the beginning of the next cycle.

In the second model, we consider a system consists of a central rework plant, which is
capable of handling imperfect quality items for rework only, and $n$ local production plants.
The behaviour of the local production plants is the same of the previous model expect they can not handle rework process.
After the production processes in the local plants, the detective items are aggregated and shipped to the central plant for rework.
The imperfect quality items are then recovered and the recovered items are for satisfying the demands at the central plant.
At the end of the cycle if there is stock at the central plant, they will be sold as one lot at a lower price.
The remainder of this paper is organized as follows. 
In Section 2, we give a literature review and the motivation of this study.
In Section 3, we consider a model for the a single production plant. 
In Section 4, we consider an aggregated model for a central rework plant and $n$ local production plants. 
We then give two numerical examples in Section 5. Finally, concluding remarks are given in Section 6 to conclude the paper.

\section{Literature Review}
Economic production quantity (EPQ) is one of the main research topic in production and inventory management.
By using EPQ model, optimal quantity of items produced can be obtained.
Classical EPQ model was developed under various assumptions.
Since then, researchers have extended the model by relaxing one or more of its assumptions.

It was assumed that the items produced is of perfect quality in the classical model.
However, imperfect quality items may be produced in reality.
\citet{Salameh} proposed an EPQ model with imperfect quality products.
The defective items are screened out and sold as a single batch at a lower price. 
\citet{Wee2007} extended the model by considering random defective rate.
\citet{Jaber} assumed the percentage defective per lot reduces according to a learning curve.
\citet{Chang2004} applied fuzzy sets theory for modeling defective rate and demands.
\citet{Rezael} considered a supply chain with multiple products and multiple suppliers.
Received items from suppliers were not of perfect quality and decision was made by using genetic algorithm.
\citet{Chung2009} proposed an inventory model with two warehouses, where one of them was rented.
\citet{Yassine} considered disaggregating the shipments of imperfect quality items in a single production run and aggregating the shipments
of imperfect items over multiple production runs.
A review on EPQ models for imperfect quality items can be found in \citep{Khan}.

How to handle imperfect quality items is another important issue.
One possible way is to perform rework and make them become serviceable.
\citet{Chan} provided a single EPQ model which considers lower pricing, rework and reject situations.
The items produced were classified as good items, good items after rework, imperfect quality items and rejected items. 
They assumed that the quality of an item was quantifiable and had a normal distribution.
\citet{Chiu1} assumed that the rework process was imperfect with random scrap rate and shortage was not allowed.
\citet{Jamal} proposed a model to obtain the optimal batch quality in a single-stage production system. 
Rework was done under two operational policies such that total system cost was minimized. 
\citet{Chiu2} determined the optimal run time for an EPQ model with scrap, rework, and stochastic machine breakdowns.
\citet{Buscher} considered a two-stage manufacturing system in which production and rework activities were carried out.
The economic production and rework quantity and the corresponding batch sizes were determined.
\citet{Taleizadeh} studied two joint production systems in a form of multiproduct single machine with and without rework.
\citet{Yoo} proposed an EPQ model that incorporated both imperfect production quality and two-way imperfect inspection.

Goods are considered as deteriorating items because their values go down with time.
Products such as electronic products, fashion clothing, food and chemical are common examples.
\citet{Teng2005} proposed an EPQ model for deteriorating items with the demand rate depended on the selling price of the products and the stock level.
\citet{Lin} considered the economic lot scheduling problem (ELSP) for deteriorating items. 
The problem was to schedule multiple products to be manufactured on a single machine repetitively over an infinite planning horizon.
\citet{Liao} developed a production model with finite production rate and considered the effect of deterioration and permissible delay in payments.
\citet{Chung2011} considered short life-cycle deteriorating items with green product design.
\citet{Widyadana} proposed an EPQ model for deteriorating items with rework, which was preformed after $m$ production setups.

One key assumption of classical EPQ model is that no shortage is allowed.
Shortage may be handled in two ways: backorders and lost sales.
\citet{Wee2006} proposed an integrated model for deteriorating items in which shortages were completely backordered.
A periodic delivery policy for a vendor and a production-inventory model for a buyer were established.
\citet{Wee2007} considered permissible shortage backordering and the effect of varying backordering cost values.
Later, \citet{Chang2010} used the renewal-reward theorem to derive the exact closed-form solutions 
of the optimal lot size, backordering quantity and maximum expected net profit per unit time.
\citet{Cardenas-Barron} extended the models in \citep{Jamal} by considering planned backorders.
Partial backordering is considered if lost sales is allowed.
\citet{Mak} proposed a optimal production-inventory policy for an inventory system with partial backordering.
\citet{Pentico} investigated the model in \citep{Mak} and redeveloped it for the EPQ with partial backordering using simpler expressions.
\citet{Wee1993} proposed an economic production policy for deteriorating items with partial
backordering using iterative method.
\citet{Giri} considered an EPQ model with increasing demand rate and adjustable production rate while shortage are partial backlogged.
\citet{Teng2007} gave a comparison between two pricing and lot-sizing models with partial backlogging and deteriorated items.

We notice that not much studies considered a model with imperfect quality and deteriorating items, rework and shortage.
On the other hand, the effect of selling deteriorated items to customers has not been addressed fully also.
In this paper, we aim at providing analytic results to address to above issues. 

\section{The Basic Models}
In this section, we consider a single production plant system which is also capable of handling imperfect quality items produced during the production process.
The following notations are used throughout the paper.\\[2mm]
\begin{tabular}{ll}
$p$ & production rate (unit/unit time)\\ 
$\al$ & percentage of good quality items produced\\
$\la$ & demand rate (unit/unit time)\\
$\th$ & percentage of items deteriorated per unit time \\ 
$\ga$ & percentage of deteriorated items screened out from the inventory\\
$p_r$ & rework process rate (unit/unit time)\\
$\al_r$ & percentage of imperfect quality items recovered\\
$\bt$ & percentage of customers who accept backlogging\\
$K$ & setup cost for a cycle (\$)\\
$c$ & deterioration cost (\$/unit)\\ 
$c_d$ & penalty cost of selling deteriorated items to customers (\$/unit)\\
$c_p$ & cost of unrecoverable imperfect quality items (\$/unit)\\
$c_s$ & shortage cost (\$/unit/unit time)\\
$c_u$ & unsatisfied demands penalty cost (\$/unit)\\ 
$h_s$ & holding cost of serviceable items (\$/unit/unit time)\\ 
$h_r$ & holding cost of imperfect quality items (\$/unit/unit time)\\[2mm]
\end{tabular}

We have the following assumptions in the development of the model.
\begin{enumerate}
\item The rates $p$, $\la$, $p_r$ and the percentages $\al$, $\ga$, $\al_r$, $\bt$ are known constants. 
\item Only serviceable items deteriorate with constant rate $\th$.
\item Shortages are allowed and are backlogged.
\item Backlogged demands are made up at the beginning of the cycle.
\item Deteriorated items and unrecoverable imperfect quality items are disposed.
\item Recovered imperfect quality items are considered as good quality items.
\end{enumerate}

The behaviour of the inventory level of serviceable items at any time during a given cycle is
illustrated in Figure 1 and the inventory level of imperfect quality items at any time during a given cycle is
illustrated in Figure 2.

\begin{figure} 
\centering 
\begin{tikzpicture} [ domain=0:15,xscale=1,yscale=.8] 
\draw [thick,<->] (0,7) -- (0,-4) -- (0,0) -- (14,0) node [below] {\; time};
\draw[thick, domain=0:1] plot (\x, {2*\x-2}); 
\draw[thick, domain=1:5] plot (\x, {(2/0.2)*(1-exp(-.2*(\x-1)))}); 
\draw[thick, domain=5:7] plot (\x, {(2/0.2*(1-exp(-.2*(5-1)))-1.5/0.2)*exp(-.2*(\x-5))+1.5/0.2}); 
\draw[thick, domain=7:11] plot (\x, {(((2/0.2*(1-exp(-.2*(5-1)))-1.5/0.2)*exp(-.2*(7-5))+1.5/0.2)+1/0.2)*exp(-.2*(\x-7))-1/0.2}); 
\draw[thick, domain=11:13.5] plot (\x, {-0.8*(\x-11)}); 
\draw[thick, dashed, domain=11:13.5] plot (\x, {-1.2*(\x-11)}); 
\draw[thick, dashed] (5,0) -- (5,5.5) -- (0,5.5) node [left] {$I_s$};
\draw[thick, dashed] (7,0) -- (7,6.163) -- (0,6.163) node [left] {$I_m$};
\node [left] at (0,-2) {$-I_b$};
\node [below] at (.4,0) {$T_1$};
\node [below] at (3,0) {$T_2$};
\node [below] at (6,0) {$T_3$};
\node [below] at (9,0) {$T_4$};
\node [below] at (12.5,0) {$T_5$};
\draw[thick, dashed] (13.5,0) -- (13.5, -3.2);
\draw [thick,<->] (0,-3.2) -- (13.5,-3.2);
\node [below] at (7,-3.2) {$T$};
\draw [thick,<->, dashed] (13.7,-3) -- (13.7,-2);
\node [align=left, right] at (13.7,-2.5) {lost\\sales};
\end{tikzpicture} 
\caption{Inventory level of serviceable items.}
\end{figure}
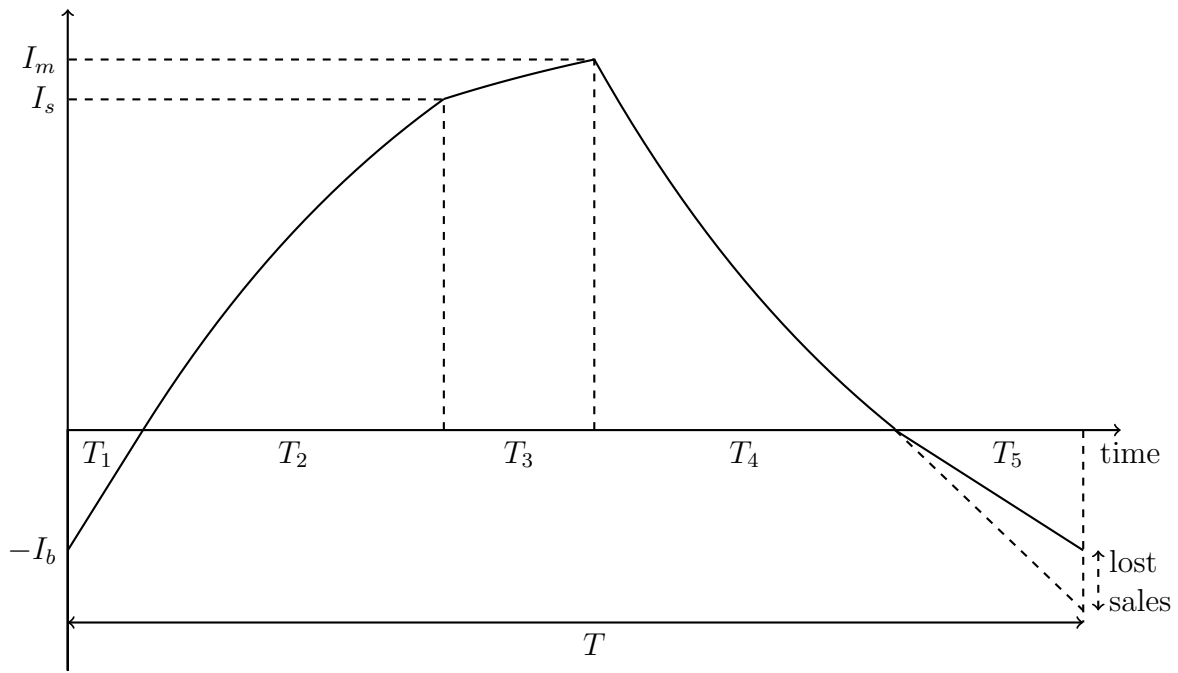 

\begin{figure} 
\centering 
\begin{tikzpicture}  [ domain=0:15,xscale=1,yscale=.8] 
\draw [thick,<->] (0,5) -- (0,0) -- (14,0) node [below] {\; time};
\draw[thick] (0,0) -- (5,4) -- (7,0);
\draw[thick, dashed] (5,0) -- (5,4) -- (0,4) node [left] {$\ \ I_c$};
\node [below] at (2.5,0) {$T_1+T_2$};
\node [below] at (6,0) {$T_3$};
\end{tikzpicture} 
\caption{Inventory level of imperfect quality items.}
\end{figure}
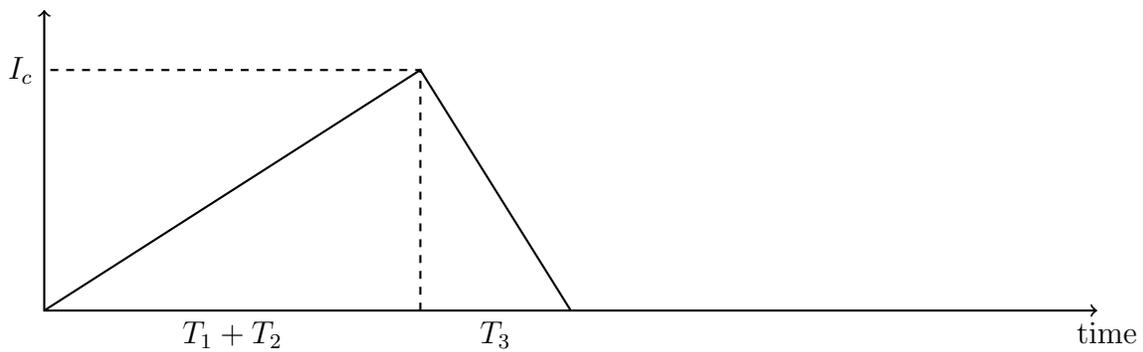 

\subsection{Model for Partial Backlogging}
The inventory level of serviceable items at time $t$ over the five periods in a cycle are determined by the following differential equations:
\begin{align}
I'_1(t_1) &= \al p -\la, \quad 0 \leq t_1 \leq T_1, \label{I1}\\
I'_2(t_2)+\gath I_2(t_2)&=\al p -\la, \quad 0 \leq t_2 \leq T_2, \label{I2}\\
I'_3(t_3)+\gath I_3(t_3)&=\al_r p_r -\la, \quad 0 \leq t_3 \leq T_3, \label{I3}\\
I'_4(t_4)+\gath I_4(t_4)&= -\la, \quad 0 \leq t_4 \leq T_4, \label{I4}\\
I'_5(t_5)&= -\bt\la, \quad 0 \leq t_5 \leq T_5. \label{I5}
\end{align}
With the boundary conditions $I_1(0)=I_5(T_5)=-I_b$, $I_1(T_1)=I_2(0)=0$, $I_2(T_2)=I_3(0)=I_s$, $I_3(T_3)=I_4(0)=I_m$, $I_4(T_4)=I_5(0)=0$, the solutions for the above differential equations are
\begin{align}
I_1(t_1)&=(\al p -\la)t_1-I_b, \quad 0 \leq t_1 \leq T_1, \label{I1t}\\
I_2(t_2)&=\Big(\dfrac{\al p -\la}{\gath}\Big)(1-\exp(-\gath t_2)), \quad 0 \leq t_2 \leq T_2, \label{I2t}\\
I_3(t_3)&=\Big(I_s-\dfrac{\al_r p_r -\la}{\gath}\Big)\exp(-\gath t_3)+\Big(\dfrac{\al_r p_r -\la}{\gath}\Big), \quad 0 \leq t_3 \leq T_3, \label{I3t}\\
I_4(t_4)&=\Big(I_m+\dfrac{\la}{\gath}\Big)\exp(-\gath t_4)- \Big(\dfrac{\la}{\gath}\Big), \quad 0 \leq t_4 \leq T_4, \label{I4t}\\
I_5(t_5)&=-\bt\la t_5, \quad 0 \leq t_5 \leq T_5.\label{I5t}
\end{align}
Hence, it can be deduced from Eq. (\ref{I2t}) that at $t_2 = T_2$,
the inventory level of serviceable items is
\begin{equation}\label{Is}
I_s=\Big(\dfrac{\al p -\la}{\gath}\Big)(1-\exp(-\gath T_2)),
\end{equation}
and from Eq. (\ref{I3t}) that at $t_3=T_3$,
the maximum inventory level serviceable items is
\begin{equation}\label{Im-Is}
I_m=\Big(\dfrac{\al_r p_r -\la}{\gath}\Big)(1-\exp(-\gath T_3))+\exp(-\gath T_3)I_s,
\end{equation}
and from Eq. (\ref{I4t}) that at $t_4=T_4$, $I_4(T_4)=0$, which gives
\begin{equation}\label{Im}
I_m=\Big(\dfrac{\la}{\gath}\Big)(\exp(\gath T_4)-1).
\end{equation}
Also, from Eqs. (\ref{I1t}) and (\ref{I5t}) the unfilled order backlog is
\begin{equation}\label{Ib}
I_b=(\al p - \la) T_1 = \bt \la T_5.
\end{equation}

For the imperfect quality items, the maximum inventory level is given by
\begin{equation}\label{Ic}
I_c=(1- \al)p (T_1+T_2) = p_r T_3.
\end{equation}

The total operating cost consists of the deterioration cost, penalty cost of selling deteriorated items to customers,
holding costs of serviceable and imperfect quality items,
setup cost for a cycle, cost of unrecoverable imperfect quality items, shortage cost and unsatisfied demands penalty cost.
Our aim is to minimize the total cost per unit time, which can be expressed as
\begin{equation}\label{TC}
\begin{split}
TC=\;&\dfrac{c}{T}((\al p-\la)T_2-I_s)+\dfrac{c}{T}((\al_r p_r-\la)T_3-(I_m-I_s))+\dfrac{c}{T}(I_m-\la T_4)\\[2mm]
&+\dfrac{1-\ga}{\ga}\Big[\dfrac{c_d}{T}((\al p-\la)T_2-I_s)+\dfrac{c_d}{T}((\al_r p_r-\la)T_3-(I_m-I_s))+\dfrac{c_d}{T}(I_m-\la T_4)\Big]\\[2mm]
&+\dfrac{h_s}{T}\Big( \dint_0^{T_2} I_2(t_2)d\, t_2 + \dint_0^{T_3} I_3(t_3)d\, t_3 + \dint_0^{T_4} I_4(t_4)d\, t_4 \Big)\\[2mm]
&+\dfrac{h_r}{T}\Big( \dfrac{(T_1+T_2+T_3)p_r T_3}{2}\Big)\\[2mm]
&+\dfrac{K}{T}+ \dfrac{c_p}{T}(1-\al_r)p_r T_3 + \dfrac{c_s}{T}\dfrac{(\al p-\la){T_1}^2}{2}+\dfrac{c_s}{T}\dfrac{\bt\la{T_5}^2}{2}+\dfrac{c_u}{T}\bt'\la T_5,
\end{split}
\end{equation}
where $\bt'=1-\bt$.

In order to simplify the expression, the relationships between $T_i$ $(i=1,\ldots,5)$ and $T$ are needed.
Firstly, $T=T_1+T_2+T_3+T_4+T_5$ together with Eq. (\ref{Ib}) give
\begin{equation}\label{T1}
T_1=\dfrac{\bt\la}{\al p-\bt'\la}(T-T_2-T_3-T_4)
\end{equation}
and
\begin{equation}\label{T5}
T_5=\dfrac{\al p -\la}{\al p-\bt'\la}(T-T_2-T_3-T_4).
\end{equation}
Then, with Eqs. (\ref{Im}) and (\ref{Ic}), Eq. (\ref{TC}) can be expressed as
\begin{equation}\label{TC1}
\begin{split}
TC=\;&\Big(\dfrac{c}{T}+ \dfrac{(1-\ga)c_d}{\ga T}\Big)((\al p-\la)T_2+(\al_r p_r-\la)T_3-\la T_4)\\[2mm]
&+\dfrac{h_s}{T}\Big[ \Big(\dfrac{\al p -\la}{(\gath)^2}\Big)(\gath T_2+\exp(-\gath T_2)-1)\\[2mm] 
&\qquad\ +\Big(\dfrac{I_s}{\gath}-\dfrac{\al_r p_r -\la}{(\gath)^2}\Big) ( 1-\exp(-\gath T_3))+  \dfrac{\al_r p_r -\la}{\gath}T_3\\[2mm] 
&\qquad\ +\Big(\dfrac{\la}{(\gath)^2}\Big) ( \exp(\gath T_4)-1-\gath T_4) \Big]\\[2mm]
&+\dfrac{h_r}{T}\Big[ \dfrac{({p_r}^2+(1-\al)p \cdot p_r) {T_3}^2}{2(1-\al)p}\Big]\\[2mm]
&+\dfrac{K}{T}+ \dfrac{c_p}{T}(1-\al_r)p_r T_3+ \dfrac{c_s}{T}\dfrac{(\al p-\la)\bt \la}{2(\al p -\bt'\la)}(T-T_2-T_3-T_4)^2\\[2mm]
&+\dfrac{c_u}{T}\dfrac{(\al p -\la)\bt'\la }{\al p-\bt'\la}(T-T_2-T_3-T_4).
\end{split}
\end{equation}

In order to find the optimal cycle length and the optimal time periods $T_i$ $(i=1,\ldots,5)$, we use the following method for getting an approximation of $TC$ for which the optimal solution can easily be obtained.
Such approximation is commonly used in modelling systems for deteriorating items, see for example \citep{Wee1993,Widyadana,Yang}.
Next, we express $T_2$ in terms of $T_3$, $T_4$ and $T$.
From Eqs. (\ref{Ic}) and (\ref{T1}), $T_2$ can be expressed as
\begin{equation}\label{T2}
T_2=\dfrac{\om}{\al p -\la}T_3+\dfrac{\bt \la}{\al p-\la}(T_4-T),
\end{equation}
where
$$
\om=\bt\la+\dfrac{(\al p - \bt'\la)p_r}{(1-\al)p}.
$$
From Eqs. (\ref{Is}), (\ref{Im-Is}) and (\ref{Im}), we have 
\begin{equation}\label{T2T3T4}
\la(\exp(\gath T_4)-1)-(\al p -\la)\exp(-\gath T_3)(1-\exp(-\gath T_2))=(\al_r p_r -\la)(1-\exp(-\gath T_3)).
\end{equation}
In what follows, we use the Taylor series approximation under the assumptions that $\gath T_2$, $\gath T_3$ and $\gath T_4$ are small.
We adopt that
$$
\exp(x)\approx 1+x+\dfrac{x^2}{2}
$$
when $x$ is small.
We simplify Eq. (\ref{T2T3T4}) to 
\begin{equation}\label{T3}
T_3=\dfrac{\la}{\al_r p_r-\la}\Big(T_4+\dfrac{\gath{T_4}^2}{2}\Big)-\dfrac{\al p-\la}{\al_r p_r-\la}T_2.
\end{equation}
The derivation can be found in the Appendix.
With Eqs. (\ref{T2}) and (\ref{T3}), the first term of $TC$ in Eq. (\ref{TC1}) then becomes
\begin{equation}\label{cd}
\dfrac{\ga c+(1-\ga)c_d}{T}\Big(\dfrac{\la\th{T_4}^2}{2}\Big).
\end{equation}

Under the assumption that $\gath T_2$ is small, $I_s$ in Eq. (\ref{Is}) becomes
$$
I_s=(\al p -\la)\Big( T_2 -\dfrac{\gath {T_2}^2}{2} \Big).
$$
Similarly, it can be shown that
the second term of $TC$ in Eq. (\ref{TC1}) is
\begin{equation}\label{hs}
\dfrac{h_s}{T}\Big[ \dfrac{(\al p-\la){T_2}^2}{2} + (\al p -\la)\Big( T_2 -\dfrac{\gath {T_2}^2}{2} \Big) \Big(T_3-\dfrac{\gath {T_3}^2}{2}\Big)+ \dfrac{(\al_r p_r -\la){T_3}^2}{2} +\dfrac{\la {T_4}^2}{2} \Big].
\end{equation}
If we further assume that $1-\gath T_2/2 \approx 1$ and $1-\gath T_3/2\approx 1$ then Eq. (\ref{hs}) can be simplified to 
\begin{equation}\label{hs1}
\dfrac{h_s}{T}\Big[ \dfrac{(\al p-\la){T_2}^2}{2} + (\al p -\la) T_2 T_3+ \dfrac{(\al_r p_r -\la){T_3}^2}{2} +\dfrac{\la {T_4}^2}{2} \Big].
\end{equation}
Hence the expression of $TC$ can be simplified as
\begin{equation}\label{TC2}
\begin{split}
TC=\;&\dfrac{\ga c+(1-\ga)c_d}{T}\Big(\dfrac{\la\th{T_4}^2}{2}\Big)\\[2mm]
&+\dfrac{h_s}{T}\Big[ \dfrac{(\al p-\la){T_2}^2}{2} + (\al p -\la) T_2 T_3+ \dfrac{(\al_r p_r -\la){T_3}^2}{2} +\dfrac{\la {T_4}^2}{2} \Big]\\[2mm]
&+\dfrac{h_r}{T}\Big[ \dfrac{({p_r}^2+(1-\al)p \cdot p_r) {T_3}^2}{2(1-\al)p}\Big]\\[2mm]
&+\dfrac{K}{T}+ \dfrac{c_p}{T}(1-\al_r)p_r T_3+ \dfrac{c_s}{T}\dfrac{(\al p-\la)\bt \la}{2(\al p -\bt'\la)}(T-T_2-T_3-T_4)^2\\[2mm]
&+\dfrac{c_u}{T}\dfrac{(\al p -\la)\bt'\la }{\al p-\bt'\la}(T-T_2-T_3-T_4).
\end{split}
\end{equation}

The optimal values of $(T_4,T)$ can be obtained by substituting $T_2$ and $T_3$ in Eqs. (\ref{T2}) and (\ref{T3}) into Eq. (\ref{TC2})
and solving $\pa TC/\pa T_4=0$ and $\pa TC/\pa T=0$ simultaneously.
Such computation can be done by any mathematical software.

\subsection{Model for Complete Backlogging}
If we assume that all the unsatisfied demand are backlogged, i.e. $\bt=1$ and $\bt'=0$, then 
the total cost per unit time can be expressed as
\begin{equation}\label{TC3}
\begin{split}
TC=\;&\dfrac{\ga c+(1-\ga)c_d}{T}\Big(\dfrac{\la\th{T_4}^2}{2}\Big)\\[2mm]
&+\dfrac{h_s}{T}\Big[ \dfrac{(\al p-\la){T_2}^2}{2} + (\al p -\la) T_2 T_3+ \dfrac{(\al_r p_r -\la){T_3}^2}{2} +\dfrac{\la {T_4}^2}{2} \Big]\\[2mm]
&+\dfrac{h_r}{T}\Big[ \dfrac{({p_r}^2+(1-\al)p \cdot p_r) {T_3}^2}{2(1-\al)p}\Big]\\[2mm]
&+\dfrac{K}{T} + \dfrac{c_p}{T}(1-\al_r)p_r T_3+ \dfrac{c_s}{T}\dfrac{(\al p-\la) \la}{2\al p }(T-T_2-T_3-T_4)^2,
\end{split}
\end{equation}
with
\begin{equation}\label{T11}
T_1=\dfrac{\la}{\al p}(T-T_2-T_3-T_4),
\end{equation}
\begin{equation}\label{T51}
T_5=\dfrac{\al p -\la}{\al p}(T-T_2-T_3-T_4)
\end{equation}
and
\begin{equation}\label{T21}
T_2=\dfrac{\om}{\al p -\la}T_3+\dfrac{\la}{\al p-\la}(T_4-T),
\end{equation}
where
$$
\om=\la+\dfrac{\al p_r}{1-\al}.
$$
For $T_3$, we use a simplified version of Eq. (\ref{T3}) under the assumption that $1+\gath T_4/2 \approx 1$: 
\begin{equation}\label{T41}
T_3=\dfrac{\la}{\al_r p_r-\la} T_4-\dfrac{\al p-\la}{\al_r p_r-\la}T_2.
\end{equation}
Together with Eq. (\ref{T21}), we have
\begin{equation}\label{T31}
T_3=\dfrac{\la}{\om+\al_r p_r -\la}T=\dfrac{(1-\al)\la}{\al p_r+(1-\al)\al_rp_r }T= \eta T.
\end{equation}

Finally, we can express $TC$ in terms of $T_4$ and $T$ only: 
\begin{equation}\label{TC4}
\begin{split}
TC(T_4,T)=\;&\dfrac{\ga c+(1-\ga)c_d}{T}\Big(\dfrac{\la\th{T_4}^2}{2}\Big)\\[2mm]
&+\dfrac{h_s}{T}\Big[ \dfrac{(\la T_4-(\al_r p_r-\la)\eta T)^2}{2(\al p-\la)} +\la \eta T_4 T - \dfrac{(\al_r p_r -\la)\eta^2 T^2}{2} +\dfrac{\la {T_4}^2}{2} \Big]\\[2mm]
&+\dfrac{h_r}{T}\Big[ \dfrac{({p_r}^2+(1-\al)p \cdot p_r) \eta^2 T^2}{2(1-\al)p}\Big]+\dfrac{K}{T}+ \dfrac{c_p}{T}(1-\al_r)p_r \eta T\\[2mm]
&+ \dfrac{c_s}{T}\dfrac{\la}{2\al p (\al p-\la)}[((1-\eta)(\al p -\la)+\eta(\al_r p_r-\la))T-\al pT_4]^2\\[2mm]
=\;&AT+BT_4+C\dfrac{{T_4}^2}{T}+\dfrac{K}{T}+D,
\end{split}
\end{equation}
where
\begin{equation}
\begin{split}
A=\;&h_s\Big[ \dfrac{(\al_r p_r-\la)^2\eta^2}{2(\al p - \la)}-\dfrac{(\al_r p_r -\la)\eta^2}{2}\Big]+h_r\Big[ \dfrac{({p_r}^2+(1-\al)p \cdot p_r) \eta^2}{2(1-\al)p}\Big]\\[2mm]
&+c_s \Big[ \dfrac{\la((1-\eta)(\al p -\la)+\eta(\al_r p_r-\la))^2}{2\al p(\al p-\la)} \Big];\\[2mm]
B=\;&h_s\Big[\la\eta-\dfrac{(\al_r p_r-\la)\la\eta}{\al p-\la}\Big]-c_s\Big[\dfrac{\la((1-\eta)(\al p -\la)+\eta(\al_r p_r-\la))}{\al p-\la}\Big];\\[2mm]
C=\;&(\ga c+(1-\ga)c_d)\Big(\dfrac{\la \th}{2}\Big)+h_s\Big[\dfrac{\la^2}{2(\al p-\la)}+\dfrac{\la}{2} \Big]+c_s\Big( \dfrac{\al p \la}{2(\al p-\la)}\Big);\\[2mm]
D=\;&c_p(1-\al_r)p_r \eta.
\end{split}
\end{equation}

The optimal pair can be obtained by the following theorem.
The proof can be found in the Appendix.
\begin{theorem}
$TC(T_4,T)$ is strictly convex. The optimal pair exists if
$$
B < 0 \quad\mbox{and}\quad 4AC>B^2
$$
and is given by
\begin{equation}\label{T4*T*}
(T_4^*,T^*)=\left( -B\sqrt{\dfrac{K}{C(4AC-B^2)}},2\sqrt{\dfrac{CK}{4AC-B^2}} \right).
\end{equation}
\end{theorem}

We may also deduce from Eqs. (\ref{Ic}) and (\ref{T31}) that the optimal production time and the economic production quantity are given respectively by
\begin{equation}\label{Tp}
T_p^*=T_1^*+T_2^*=\dfrac{\la}{\al p+(1-\al)\al_r p}T^*
\end{equation}
and
\begin{equation}\label{Q*}
Q^*=pT_p^*=\dfrac{\la}{\al +(1-\al)\al_r}T^*.
\end{equation}

\section{The Aggregated Model}
In this section, we consider a system consists of a central rework plant, which is capable of handling imperfect quality items for rework only, and $n$ local production plants.
After the production processes in the local plants, the imperfect quality items are aggregated and shipped to the central plant for rework. 
We assume that all the local production plants are the same as the one considered in Section 3 except rework is not preformed there.
Therefore, the time between two shipments of imperfect quality items is $T$.
It is natural to use $T$ as the cycle length for the central rework plant also.
Besides the notations and assumptions given in Section 3, we have more as stated in the following.

\noindent Notations:\\[2mm]
\begin{tabular}{ll}
$K_c$ & setup cost for a cycle at the central rework plant (\$)\\
$c_v$ & penalty cost of selling recovered items at the end of the
cycle (\$/unit)\\ 
$h_c$ & holding cost at the central rework plant (\$/unit/unit time)\\[2mm]
\end{tabular}

\noindent Assumptions:
\begin{enumerate}
\item Shortages were completely backlogged at the local production plants.
	\item The transportation cost for the imperfect quality items is included in $K$.
	\item The transportation time for the imperfect quality items is neglected.
	\item Rework is processed in no time and all imperfect quality items are recovered to good quality items.
\end{enumerate}

At each of the local production plant, the behaviour of the inventory level of serviceable items at any time during a given cycle is
illustrated in Figure 3 and the inventory level of imperfect quality items at any time during a given cycle is
illustrated in Figure 4.

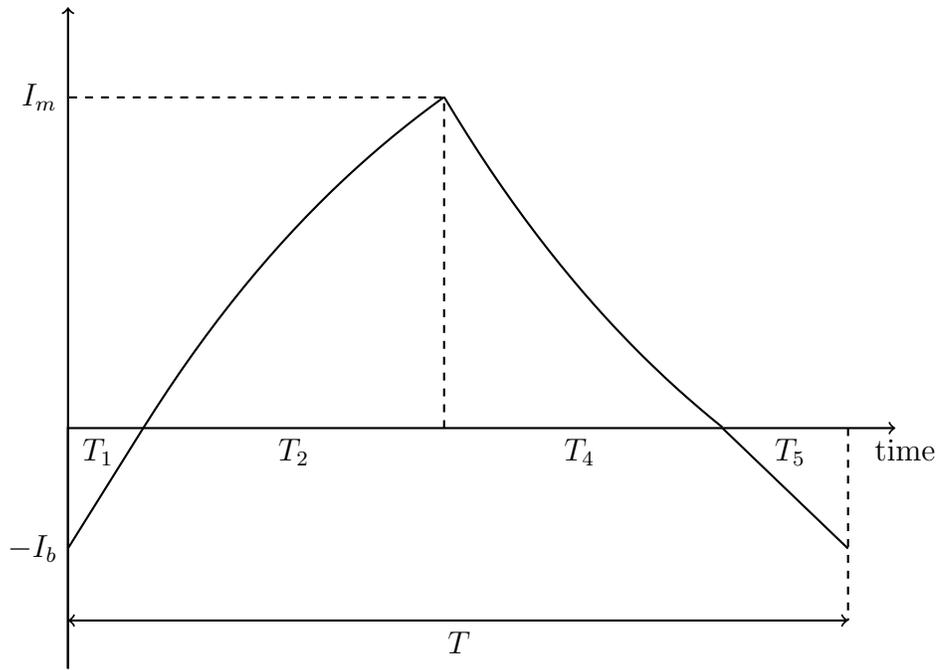
\begin{figure} 
\centering 
\begin{tikzpicture} [ domain=0:15,xscale=1,yscale=.8] 
\draw [thick,<->] (0,7) -- (0,-4) -- (0,0) -- (11,0) node [below] {\; time};
\draw[thick, domain=0:1] plot (\x, {2*\x-2}); 
\draw[thick, domain=1:5] plot (\x, {(2/0.2)*(1-exp(-.2*(\x-1)))}); 
\draw[thick, domain=5:8.7] plot (\x, {((2/0.2*(1-exp(-.2*(5-1)))+1/0.2))*exp(-.2*(\x-5))-1/0.2}); 
\draw[thick, domain=8.7:10.37] plot (\x, {-1.2*(\x-8.7)}); 
\draw[thick, dashed] (5,0) -- (5,5.5) -- (0,5.5) node [left] {$I_m$};
\node [left] at (0,-2) {$-I_b$};
\node [below] at (.4,0) {$T_1$};
\node [below] at (3,0) {$T_2$};
\node [below] at (6.8,0) {$T_4$};
\node [below] at (9.6,0) {$T_5$};
\draw[thick, dashed] (10.37,0) -- (10.37, -3.2);
\draw [thick,<->] (0,-3.2) -- (10.37,-3.2);
\node [below] at (5.2,-3.2) {$T$};
\end{tikzpicture} 
\caption{Inventory level of serviceable items at local production plants.}
\end{figure} 

\begin{figure} 
\centering 
\begin{tikzpicture}  [ domain=0:15,xscale=1,yscale=.8] 
\draw [thick,<->] (0,5) -- (0,0) -- (11,0) node [below] {\; time};
\draw[thick,thick] (0,0) -- (5,4);
\draw[thick, dashed] (5,0) -- (5,4) -- (0,4) node [left] {$\ \ I_c$};
\node [below] at (2.5,0) {$T_1+T_2$};
\end{tikzpicture} 
\caption{Inventory level of imperfect quality items at local production plants.}
\end{figure}
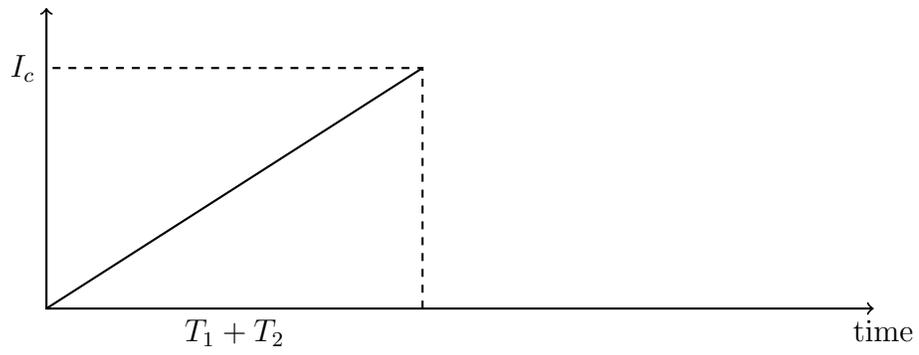 

By using the similar arguments in Section 3, we obtain
\begin{align}
I_m & = \dfrac{\al p-\la}{\gath}(1-\exp(-\gath T_2))= \dfrac{\la}{\gath}(\exp(\gath T_4)-1), \label{AIm} \\ 
I_b &= (\al p - \la)T_1 = \la T_5, \label{AIb}  \\
I_c &= (1-\al)p(T_1+T_2), \label{AIc}\\
T &= T_1+T_2+T_4+T_5.
\end{align}
Under the assumption that $\gath T_4$ is small, we get
\begin{align}
T_1 &= \dfrac{\la}{\al p}(T-T_2-T_4), \label{AT1} \\ 
T_5 &= \dfrac{\al p-\la}{\al p}(T-T_2-T_4), \label{AT5}
\end{align}
and under the assumption that $\gath T_2$ is small, $T_2$ can be approximated by 
\begin{equation}\label{AT2}
T_2 \approx \dfrac{\la}{\al p - \la} \Big(T_4+\dfrac{\gath{T_4}^2}{2}\Big).
\end{equation}
The optimization problem is to minimize the total cost per unit time, which can be expressed as
\begin{equation}\label{ATC}
	\begin{split}
TC&=n\Big[\dfrac{\ga c+(1-\ga)c_d}{T}\Big(\dfrac{\la\th{T_4}^2}{2}\Big)+\dfrac{h_s}{T} \Big( \dfrac{\la}{\al p - \la} \Big) \Big(\dfrac{\al p {T_4}^2 + \gath \la {T_4}^3}{2}\Big)\\[2mm]
&+\dfrac{h_r}{T}\dfrac{(1-\al)p}{2} \Big( \dfrac{\la T}{\al p} +\dfrac{\gath \la {T_4}^2}{2\al p} \Big)^2 \\[2mm]
&+\dfrac{K}{T}+ \dfrac{c_s}{T}\dfrac{\la}{2}\Big(\dfrac{\al p-\la}{\al p}\Big) \Big(T-\dfrac{\al p T_4+\gath \la {T_4}^2/2}{\al p - \la}\Big)^2 \Big]\\[2mm]
&+\dfrac{F(T)}{T},
\end{split}
\end{equation}
where $F(T)$ is the total cost per cycle of the central rework plant which is defined below.

By using the similar arguments in Section 3 again, the inventory level after rework at the central rework plant is given by
\begin{equation}\label{I6}
\begin{split}
I_6(t_6)&=\Big( nI_c+\dfrac{\la}{\gath} \Big)\exp(-\gath t_6)-\dfrac{\la}{\gath}, \quad 0 \leq t_6 \leq T_6,\\
&\approx\Big( nI_c+\dfrac{\la}{\gath} \Big)(1-\gath t_6)-\dfrac{\la}{\gath}\\
&=nI_c(1-\gath t_6)-\la t_6.
\end{split}
\end{equation}
We also have $I_6(T_6)=0$, which means
\begin{equation}\label{T6}
T_6=\dfrac{nI_c}{\la+nI_c\gath}.
\end{equation}

According to the sale situation, we calculate the total cost per cycle under the following two cases, which are illustrated in Figures 5 and 6 respectively:
\begin{itemize}
	\item Case I: $T_6 \geq T$, i.e. $nI_c(1-\gath T)-\la T \geq 0$.
	
	\noindent In this case, there are recovered items left in the inventory at the end of the cycle. 
	The total cost per cycle of the central rework plant is given by
	\begin{equation}\label{FT1}
	\begin{split}
	F(T) &= h_c \Big( \dint_0^T I_6(t_6)\, dt_6 \Big)+c_v(nI_c(1-\gath T)-\la T)+K_c\\
	&= h_c \Big[ nI_c T -\Big( \dfrac{\la+nI_c\gath}{2} \Big) T^2 \Big]+c_v(nI_c(1-\gath T)-\la T)+K_c.
	\end{split}
	\end{equation} 
	
	\item Case II: $T_6 < T$, i.e. $nI_c(1-\gath T)-\la T < 0$.
	
	\noindent In this case, the inventory level of recovered items become zero before the end of the cycle.
	The total cost per cycle of the central rework plant is given by
	\begin{equation}\label{FT2}
	\begin{split}
	F(T) &= h_c \Big( \dint_0^{T_6} I_6(t_6)\, dt_6 \Big)+c_u\la(T-T_6)+K_c\\
	&= h_c \Big[ nI_c T_6 -\Big( \dfrac{\la+nI_c\gath}{2} \Big) {T_6}^2 \Big]+c_u\la(T-T_6)+K_c\\
	&= h_c \Big[ \dfrac{(nI_c)^2}{2(\la+nI_c\gath)} \Big]+c_u\la(T-T_6)+K_c.
	\end{split}
	\end{equation} 
\end{itemize}

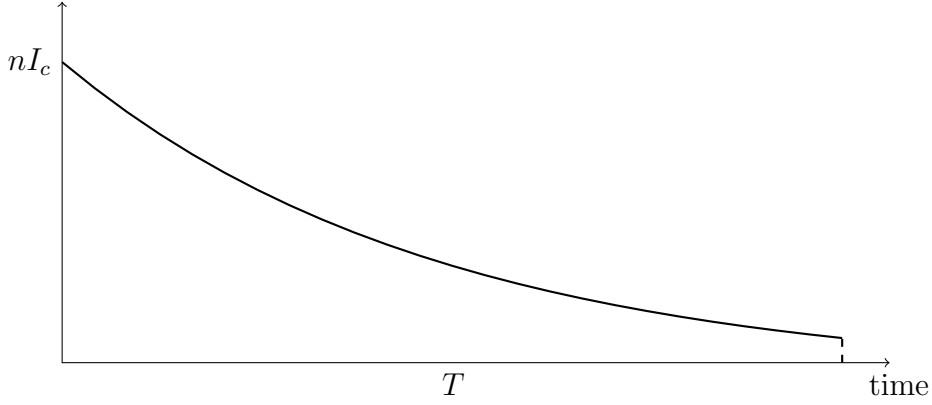
\begin{figure} 
\centering 
\begin{tikzpicture} [ domain=0:15,xscale=1,yscale=.8] 
\draw [<->] (0,6) -- (0,0) -- (11,0) node [below] {\; time};
\draw[thick, domain=0:10.37] plot (\x, {(5+0.05/0.2)*exp(-.2*\x)-0.05/0.2}); 
\draw[thick, dashed] (10.37,0) -- (10.37, 0.41);
\node [below] at (5.2,0) {$T$};
\node [left] at (0,5) {$nI_c$};
\end{tikzpicture} 
\caption{Inventory level of recovered items at central rework plant (Case I).}
\end{figure} 

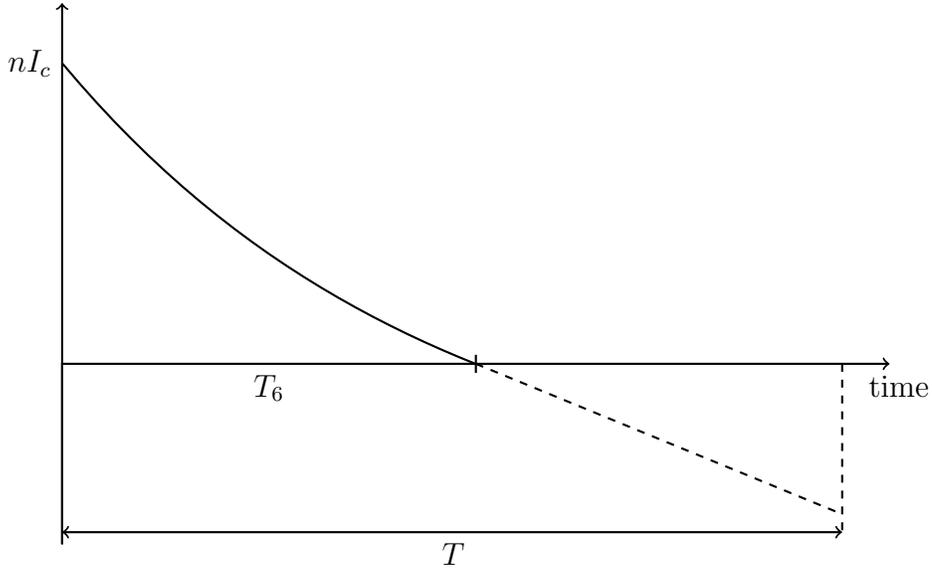
\begin{figure} 
\centering 
\begin{tikzpicture} [ domain=0:15,xscale=1,yscale=.8] 
\draw [thick,<->] (0,6) -- (0,-3) -- (0,0) -- (11,0) node [below] {\; time};
\draw[thick, domain=0:5.5] plot (\x, {(5+0.5/0.2)*exp(-.2*\x)-0.5/0.2}); 
\draw[thick, dashed] (5.5,0) -- (10.37, -2.5);
\draw[thick, dashed] (10.37,0) -- (10.37, -2.8);
\draw[thick] (5.5,0.15) -- (5.5, -0.15);
\node [below] at (2.75,0) {$T_6$};
\node [left] at (0,5) {$nI_c$};
\draw [thick,<->] (0,-2.8) -- (10.37,-2.8);
\node [below] at (5.2,-2.8) {$T$};
\end{tikzpicture} 
\caption{Inventory level of recovered items at central rework plant (Case II).}
\end{figure}

The maximum inventory level of the recovered items is $nI_c$, which can be found by Eqs. (\ref{AIc}), (\ref{AT1}) and (\ref{AT2}):
\begin{equation}\label{nIc}
nI_c=\dfrac{n\la(1-\al)}{\al}\Big( T+\dfrac{\gath{T_4}^2}{2}\Big).
\end{equation}

Now the total cost per unit time can be expressed in terms of $T_4$ and $T$ only.
However, the expression is complicated and analysis can not be done easily.
Using the fact that $\gath T_4$ and $\gath T_6$ are relatively small, we may express Eq. (\ref{ATC}) in a more simple way.
Firstly, Eqs. (\ref{AT2}), (\ref{T6}) and (\ref{nIc}) become
\begin{equation}\label{ATa}
T_2 \approx \dfrac{\la}{\al p - \la} T_4, \quad T_6 \approx \dfrac{nI_c}{\la}, \quad nI_c \approx \dfrac{n\la(1-\al)}{\al} T
\end{equation}
and
\begin{equation}\label{Con}
nI_c(1-\gath T)-\la T \geq 0 \Rightarrow \dfrac{1}{\gath}\Big( 1-\dfrac{\al}{n(1-\al)} \Big)\geq T.   
\end{equation}
Secondly, Eq. (\ref{ATC}) becomes
\begin{equation}\label{ATC1}
	\begin{split}
	TC&= n\Big[\dfrac{\ga c+(1-\ga)c_d}{T}\Big(\dfrac{\la\th{T_4}^2}{2}\Big)+\dfrac{h_s}{T} \Big( \dfrac{\la}{\al p - \la} \Big) \Big(\dfrac{\al p {T_4}^2}{2}\Big)\\[2mm]
	&\quad +\dfrac{h_r}{T}\dfrac{(1-\al)p}{2} \Big( \dfrac{\la T}{\al p} \Big)^2 +\dfrac{K}{T}+ \dfrac{c_s}{T}\dfrac{\la}{2}\Big(\dfrac{\al p-\la}{\al p}\Big) \Big(T-\dfrac{\al p T_4}{\al p - \la}\Big)^2 	\Big]\\[2mm]
	&\quad +\dfrac{F(T)}{T}.
	\end{split}
\end{equation}
Thirdly, $F(T)$ becomes
\begin{equation}\label{FT}
\begin{split}
&\qquad F(T)\\
&=\left\{\begin{array}{l}
	h_c \Big[ \Big(\dfrac{n\la(1-\al)}{\al}  -\dfrac{\la}{2} \Big) T^2 \Big]+c_v\Big(\dfrac{n\la(1-\al)}{\al}(T-\gath T^2)-\la T\Big)+K_c,\qquad \qquad\\[2mm] 
	\hfill\mbox{if }\dfrac{1}{\gath}\Big(1-\dfrac{\al}{n(1-\al)}\Big)\geq T \mbox{ (Case I)};\\[8mm]
	h_c \Big[ \dfrac{\la(n(1-\al)T/\al)^2}{2} \Big]+c_u\Big(\la T-\dfrac{n\la(1-\al)}{\al}T\Big)+K_c,\\[2mm] 
	\hfill\mbox{if }\dfrac{1}{\gath}\Big(1-\dfrac{\al}{n(1-\al)}\Big)< T \mbox{ (Case II)}.
	\end{array}\right.
	\end{split}
\end{equation}
Since $F(T)$ is a piecewise function, we tackle the optimization problem by considering two separate objective functions, which are corresponding to Case I and Case II respectively:
\begin{equation}\label{ATC2}
	\begin{split}
	TC_1&= n\Big[\dfrac{\ga c+(1-\ga)c_d}{T}\Big(\dfrac{\la\th{T_4}^2}{2}\Big)+\dfrac{h_s}{T} \Big( \dfrac{\la}{\al p - \la} \Big) \Big(\dfrac{\al p {T_4}^2}{2}\Big)\\[2mm]
	&\quad +\dfrac{h_r}{T}\dfrac{(1-\al)p}{2} \Big( \dfrac{\la T}{\al p} \Big)^2 +\dfrac{K}{T}+ \dfrac{c_s}{T}\dfrac{\la}{2}\Big(\dfrac{\al p-\la}{\al p}\Big) \Big(T-\dfrac{\al p T_4}{\al p - \la}\Big)^2 	\Big]\\[2mm]
	&\quad +\dfrac{h_c}{T} \Big[ \Big(\dfrac{n\la(1-\al)}{\al}  -\dfrac{\la}{2} \Big) T^2 \Big]+\dfrac{c_v}{T}\Big(\dfrac{n\la(1-\al)}{\al}(T-\gath T^2)-\la T\Big)+\dfrac{K_c}{T}\\[2mm]
	&= A_1T+BT_4+C\dfrac{{T_4}^2}{T}+\dfrac{nK+K_c}{T}+D_1;
	\end{split}
\end{equation}
\begin{equation}\label{ATC3}
	\begin{split}
	TC_2&= n\Big[\dfrac{\ga c+(1-\ga)c_d}{T}\Big(\dfrac{\la\th{T_4}^2}{2}\Big)+\dfrac{h_s}{T} \Big( \dfrac{\la}{\al p - \la} \Big) \Big(\dfrac{\al p {T_4}^2}{2}\Big)\\[2mm]
	&\quad +\dfrac{h_r}{T}\dfrac{(1-\al)p}{2} \Big( \dfrac{\la T}{\al p} \Big)^2 +\dfrac{K}{T}+ \dfrac{c_s}{T}\dfrac{\la}{2}\Big(\dfrac{\al p-\la}{\al p}\Big) \Big(T-\dfrac{\al p T_4}{\al p - \la}\Big)^2 	\Big]\\[2mm]
	&\quad +\dfrac{h_c}{T} \Big[ \dfrac{\la(n(1-\al)T/\al)^2}{2} \Big]-\dfrac{c_u}{T}\Big(\dfrac{n\la(1-\al)}{\al}(T-\gath T^2)-\la T\Big)+\dfrac{K_c}{T}\\[2mm]
	&= A_2T+BT_4+C\dfrac{{T_4}^2}{T}+\dfrac{nK+K_c}{T}+D_2,
	\end{split}
\end{equation}
where
\begin{equation}
	\begin{split}
	A_1&= h_r\Big(\dfrac{n(1-\al)\la^2}{2\al^2 p} \Big)+c_s\Big(\dfrac{n(\al p-\la)\la}{2\al p}\Big) 
	+h_c \Big(\dfrac{n\la(1-\al)}{\al}  -\dfrac{\la}{2} \Big)
	-c_v\Big(\dfrac{n\la\gath(1-\al)}{\al}\Big);\\
	A_2&= h_r\Big(\dfrac{n(1-\al)\la^2}{2\al^2 p} \Big)+c_s\Big(\dfrac{n(\al p-\la)\la}{2\al p}\Big) 
	+h_c \Big( \dfrac{n^2\la(1-\al)^2}{2\al^2} \Big);\\
	B&=-c_s n\la;\\
	C&=(\ga c+(1-\ga)c_d)\Big(\dfrac{n\la\th}{2}\Big)+h_s \Big( \dfrac{n\la \al p}{2(\al p - \la)} \Big) 
	+c_s \Big( \dfrac{n \la \al p}{2(\al p - \la)} \Big);\\
	D_1&=c_v\Big(\dfrac{n\la(1-\al)}{\al}-\la \Big);\\
	D_2&=c_u\Big(\la-\dfrac{n\la(1-\al)}{\al} \Big).
	\end{split}
\end{equation}

By Theorem 1, it is not difficult to see that both functions are strictly convex.
The optimal pairs which minimize $TC_1$ and $TC_2$, denoted by $(T_4^{1*},T^{1*})$ and $(T_4^{2*},T^{2*})$ respectively, are
\begin{equation}\label{T4*T*1}
(T_4^{i*},T^{i*})=\left( -B\sqrt{\dfrac{nK+K_c}{C(4A_iC-B^2)}},2\sqrt{\dfrac{C(nK+K_c)}{4A_iC-B^2}} \right), \quad i=1,2.
\end{equation}

Hence, the following procedure of finding the optimal pair is developed:
\begin{itemize}
	\item Step 1. Solve two optimization problems with objective functions $TC_1$ and $TC_2$.
  \item Step 2. If $(1/\gath)[1-\al/(n(1-\al))] \leq 0$, then $(T_4^{2*},T^{2*})$ is the optimal solution to the problem.
	\item Step 3. If $(1/\gath)[1-\al/(n(1-\al))] < T^{1*}$, which means $T^{1*}$ does not satisfy the condition of Case I, then the optimal pairs must lie on the boundary of the feasible region. Hence we replace $T^{1*}$ by $(1/\gath)[1-\al/(n(1-\al))]$  and $T_4^{1*}$ by $(-B/2C)(1/\gath)[1-\al/(n(1-\al))]$.
	\item Step 4. Similarly, if $(1/\gath)[1-\al/(n(1-\al))] \geq T^{2*}$, which means $T^{2*}$ does not satisfy the condition of Case II,  
	then replace $T^{2*}$ by $(1/\gath)[1-\al/(n(1-\al))]$ and $T_4^{2*}$ by $(-B/2C)(1/\gath)[1-\al/(n(1-\al))]$.
	\item Step 5. Compare $TC_1(T_4^{1*},T^{1*})$ with $TC_2(T_4^{2*},T^{2*})$. If $TC_1(T_4^{1*},T^{1*})$ is greater than $TC_2(T_4^{2*},T^{2*})$ then the optimal solution to the problem is
	$(T_4^*,T^*)=(T_4^{1*},T^{1*})$. Otherwise, the optimal solution to the problem is $(T_4^*,T^*)=(T_4^{2*},T^{2*})$.
\end{itemize}

We may also deduce from Eqs. (\ref{AT1}) and (\ref{AT2}) that the optimal production time and the economic production quantity of the local production plants are given respectively by
\begin{equation}\label{ATp}
T_p^*=T_1^*+T_2^*=\dfrac{\la}{\al p}\Big(T^*+\dfrac{\gath {T_4^*}^2}{2}\Big)
\end{equation}
and
\begin{equation}\label{AQ*}
Q^*=pT_p^*=\dfrac{\la}{\al}\Big(T^*+\dfrac{\gath {T_4^*}^2}{2}\Big).
\end{equation}

\section{Numerical Examples}
In this section, we provide numerical examples for the models developed in the previous sections.
We first give an example for the single production plant model in Section 3.
We assume that $p=6000,$ $\al=0.7$, $\th=0.1$, $\ga=0.6$, $\la=1000$, $p_r=4000$ and $\al_r=0.6$.
For the operating costs, we assume that $c=\$40$, $c_p=\$30$, $c_s=\$200$, $c_d=\$100$, $h_s=\$5$, $h_r=\$4$ and $K=\$300$.
We follow from Eq. (\ref{T4*T*}) that the optimal pair is given by
$$
(T_4^*,T^*)=(0.1996,0.2891).
$$
We find $T_2^*$ and $T_3^*$ by solving Eqs. (\ref{T3}) and (\ref{T21}):
$$T_2^*=0.0519, \quad T_3^*=0.0247.$$
By Eqs. (\ref{T11}) and (\ref{T51}), we have
$$
T_1^*=0.0031, \quad T_5^*=0.0098.
$$
From Eqs. (\ref{Tp}) and (\ref{Q*}) the optimal production time and the economic production quantity are given respectively by
$$
T_p^*=0.0550 \quad \mbox{and} \quad Q^* \approx 330.
$$
From Eqs. (\ref{Is}) and (\ref{Im}) the maximum inventory is 
$$
I_m \approx 201
$$
while the maximum inventory during production period is
$$
I_s \approx 166.
$$
From Eqs. (\ref{Ib}) and (\ref{Ic}) the amount of shortages per cycle is
$$
I_b \approx 10
$$
while the maximum inventory level for imperfect quality items is
$$
I_c \approx 99.
$$
The corresponding optimal total cost per unit time can be found by using Eq. (\ref{TC1}):
$$
TC^*=\$5837.6.
$$

We remark that although $(T_4^*,T^*)=(0.1996,0.2891)$ is the optimal pair to the problem with objective function (\ref{TC4}) rather than the original problem 
with objective function (\ref{TC1}), it still gives a reasonable approximation to the original optimal pair.
In fact, using numerical method, the original optimal pair is $(T_4^*,T^*)=(0.2449,0.3493)$ with $TC^*=\$5800.7$, which means the approximation error is less than $1\%$.


We next give an example for the aggregated plant model to demonstrate the use of the procedure developed in Section 4.
We assume the same parameters as the previous example with $n=5$, $K_c=\$250$, $c_v=\$10$, $h_c=\$3$.
We follow from Eq. (\ref{T4*T*1}) that the optimal pairs are given by
$$
(T_4^{1*},T^{1*})=(0.2469,0.3337), \quad (T_4^{2*},T^{2*})=(0.2248,0.3038).
$$
We next check whether $(1/\gath)[1-\al/(n(1-\al))] < T^{1*}$ and $(1/\gath)[1-\al/(n(1-\al))] \geq T^{2*}$.
We have
$$
\dfrac{1}{\gath}\Big( 1-\dfrac{\al}{n(1-\al)} \Big)=8.8889,
$$
which means we need to replace $(T_4^{2*},T^{2*})$ by
$$
(T_4^{2*},T^{2*})=(6.5760,8.8889).
$$
Finally we compare $TC_1(T_4^{1*},T^{1*})=\$21917$ with $TC_2(T_4^{2*},T^{2*})=\$145868$,
therefore we choose
$$
(T_4^*,T^*)=(0.2469,0.3337).
$$
From Eq. (\ref{nIc}) the maximum inventory level for recovered items is
$$
nI_c \approx 719.
$$
From Eqs. (\ref{ATp}) and (\ref{AQ*}) optimal production time and the economic production quantity of the local production plants are given respectively by
$$
T_p^*=0.0799\quad \mbox{and} \quad Q^* \approx 479.
$$

\section{Concluding Remarks}
In this paper, two economic production quantity (EPQ) models were proposed for deteriorating items
with rework process.
We first developed a model for a single production plant system which is also capable of processing rework.
We then developed an aggregated model for a system consists of $n$ local production plants and one central rework plant.
The approximated optimal cycle length and the EPQ were obtained.
Numerical examples were also provided to illustrate the solution procedure.
In this paper, the demand rate was assumed to be a constant.
An extension to this paper can be done by considering time-dependent demand rate or stochastic demand rate.

\section*{Appendix}

\noindent{\bf Derivation of Eq. (\ref{T3}).}
We have (Eq. (\ref{T2T3T4}))
$$
\la(\exp(\gath T_4)-1)-(\al p -\la)\exp(-\gath T_3)(1-\exp(-\gath T_2))=(\al_r p_r -\la)(1-\exp(-\gath T_3)).
$$
By using the exponential series approximation, we have 
\begin{align*}
&\la\Big(\gath T_4 +\dfrac{(\gath T_4)^2}{2} \Big)-(\al p - \la)\Big(1-\gath T_3 +\dfrac{(\gath T_3)^2}{2} \Big)
\Big(\gath T_2 -\dfrac{(\gath T_2)^2}{2} \Big)\\
=&(\al_r p_r -\la) \Big(\gath T_3 -\dfrac{(\gath T_3)^2}{2} \Big)\\
&\la\Big( T_4 +\dfrac{\gath {T_4}^2}{2} \Big)-(\al p - \la)\Big(1-\gath T_3 +\dfrac{(\gath T_3)^2}{2} \Big)
\Big( T_2 -\dfrac{\gath{T_2}^2}{2} \Big)\\
=&(\al_r p_r -\la) \Big( T_3 -\dfrac{\gath {T_3}^2}{2} \Big).
\end{align*}
Under the assumptions that $\gath T_2$ and $\gath T_3$ are small, we have 
$$
\la\Big( T_4 +\dfrac{\gath {T_4}^2}{2} \Big)-(\al p - \la)T_2=(\al_r p_r -\la) T_3,
$$
which means
$$
T_3=\dfrac{\la}{\al_r p_r-\la}\Big(T_4+\dfrac{\gath{T_4}^2}{2}\Big)-\dfrac{\al p-\la}{\al_r p_r-\la}T_2,
$$
as desired. \qed

\noindent{\bf Proof of Theorem 1.}
Differentiating $TC(T_4,T)=AT+BT_4+C\dfrac{{T_4}^2}{T}+\dfrac{K}{T}+D$ with respect to $T_4$ and $T$, we have
$$
\begin{array}{rcl}
\dfrac{\pa TC}{\pa T_4} &=& B+\dfrac{2CT_4}{T};\\[2mm]
\dfrac{\pa TC}{\pa T} &=& A-\dfrac{K}{T^2}-\dfrac{C{T_4}^2}{T^2}.\\[2mm]
\end{array}
$$
Solving $\pa TC/\pa T_4=0$ and $\pa TC/\pa T=0$ simultaneously, we get
$$
(T_4^*,T^*)=\left( -B\sqrt{\dfrac{K}{C(4AC-B^2)}},2\sqrt{\dfrac{CK}{4AC-B^2}} \right)
$$
as desired.
Since we must have $T_4^*>0$ and $T^*>0$, therefore
$$
B < 0 \quad\mbox{and}\quad 4AC>B^2.
$$

The Hessian matrix is given by
$$
H=\left(\begin{array}{cc}
\dfrac{2C}{T}&-\dfrac{2CT_4}{T}\\[2mm]
-\dfrac{2CT_4}{T}&\dfrac{2K}{T^3}+\dfrac{2CT_4^2}{T^3}
\end{array}
\right).
$$
It is easy to see that since $2C/T>0$ and $\det(H)=4CK/T^4>0$, $H$ is positive definite and hence 
$TC(T_4,T)$ is strictly convex.
\qed

\end{document}